\documentclass[a4paper,11pt,reqno]{amsart}

\usepackage[full]{textcomp}

\usepackage{charter}

\usepackage{graphicx}
\usepackage{color}
\definecolor{rblue}{RGB}{39,64,139}
\usepackage{hyperref}
\hypersetup{pdftex,naturalnames=false,%
	draft=false,
	bookmarks=true,bookmarksopenlevel=1,%
	allcolors=rblue,
	allbordercolors={.15 .25 .55},
	colorlinks=true,
	}
\usepackage{amssymb}
\usepackage[initials]{amsrefs}
\usepackage{mathtools}
\usepackage[cmtip]{xypic}
\title[Lie-algebraic approach to the local index theorem]{A Lie-algebraic approach to the local index theorem on compact homogeneous spaces}
\author{Seunghun Hong}
\address{Northwestern College, 101 7th St SW, Orange City, 51041 IA, U.S.A.}
\email{seunghun.hong@nwciowa.edu}
\urladdr{diracoperat.org}


\theoremstyle{plain}
\newtheorem{thm}{Theorem}[section]
\newtheorem{lem}[thm]{Lemma}
\newtheorem{prop}[thm]{Proposition}

\theoremstyle{definition}

\newtheorem{notation}[thm]{Notation}
\theoremstyle{remark}
\newtheorem*{rmk}{Remark}

\newcommand{\KD}{\mathcal{D}}
\newcommand{\C}{\mathbb{C}}
\newcommand{\R}{\mathbb{R}}
\newcommand{\Z}{\mathbb{Z}}

\newcommand{\bas}{\mathrm{bas}}
\newcommand{\topf}{\mathrm{top}}
\newcommand{\spin}{\mathfrak{spin}}
\newcommand{\so}{\mathfrak{so}}
\newcommand{\qW}{\mathcal{W}}
\newcommand{\Id}{\mathbf{1}}
\DeclareMathOperator{\Spin}{Spin}

\DeclareMathOperator{\SU}{SU}
\DeclareMathOperator{\Str}{Str}
\DeclareMathOperator{\tr}{tr}
\newcommand{\vol}{\mathrm{vol}}
\DeclareMathOperator{\End}{End}
\DeclareMathOperator{\Aut}{Aut}
\DeclareMathOperator{\Ind}{Ind}
\DeclareMathOperator{\ad}{ad}
\DeclareMathOperator{\Cl}{Cl}
\DeclareMathOperator{\Ad}{Ad}

\DeclareMathOperator{\diag}{diag}
\DeclareMathOperator{\Duf}{\mathcal{Q}}

\DeclareMathOperator{\ch}{ch}
\DeclareMathOperator{\CW}{CW}
\DeclareMathOperator{\pr}{pr}
\DeclareMathOperator{\ev}{ev}
\DeclareMathOperator{\Taylor}{Taylor}

\newcommand{\ldef}{\coloneqq}

\begin{document}
\begin{abstract}
Using a K-theory point of view, \mbox{Bott} related the \mbox{Atiyah}-\mbox{Singer} index theorem for elliptic operators on compact homogeneous spaces to the Weyl character formula. This article explains how to prove the local index theorem for compact homogenous spaces using Lie algebra methods. The method follows in outline the proof of the local index theorem due to \mbox{Berline} and \mbox{Vergne}. But the use of \mbox{Kostant's} cubic Dirac operator in place of the Riemannian Dirac operator leads to substantial simplifications. An important role is also played by the quantum Weil algebra of \mbox{Alekseev} and \mbox{Meinrenken}. 
\end{abstract}
\maketitle

\section{Introduction}

Soon after \mbox{Atiyah} and \mbox{Singer}~\cite{atiyahsingerann} proved the index theorem, \mbox{Bott}~\cite{bott} examined it over homogeneous spaces $G/K$ where $G$ is a compact connected Lie group and $K$ is a closed connected subgroup. Using \mbox{Weyl's} theory---integral formula, character formula, and so on---\mbox{Bott} was able to verify the index theorem by-passing analytic or topological arguments considerably. Expectation of such simplification for the \emph{local} index theorem is the motive of this article. 

It seems that the well-known proofs for the local index theorem (such as the ones found in Atiyah, Bott, and Patodi~\cite{abp}, Getzler~\cite{getzler}, Bismut~\cite{bismut}, or Berline and Vergne~\cite{berlinevergne}) in themselves do not become simpler even if we restrict the manifolds under consideration to compact homogeneous spaces. A common feature of the aforementioned proofs is that they use the Riemannian Dirac operator. What we shall demonstrate is that the use of \mbox{Kostant's} \emph{cubic Dirac operator}~\cite{kostant} in place of the Riemannian Dirac operator leads to substantial simplifications; \mbox{Alekseev} and \mbox{Meinrenken's} \emph{quantum Weil algebra}~\cites{alekmein,alekmein2} also plays an important role.

This is certainly not the first instance where the utility of the cubic Dirac operator is found. Its usefulness has already been demonstrated in the realm of the representation theory of complex semisimple Lie algebras. We do not list all the works in that direction but refer to the treatise by Huang and \mbox{Pand{\v{z}}i{\'c}}~\cite{huang}. On a more differential geometric side, the advantage of the cubic Dirac operator for index theoretic purposes has been indicated by the works of Slebarski~\cites{slebarski1,slebarski2,slebarski3} and Goette~\cite{goette}. More recently, Freed, Hopkins, and Teleman~\cite{fht2} used a family of cubic Dirac operators in the context of loop groups to prove the ``Thom isomorphism'' in that setting.

We shall outline our approach in the next few paragraphs. But first let us recapitulate the local index theorem. Let $\KD$ be a Dirac operator on a vector bundle $S$ over a closed even-dimensional manifold $M$. Owing to the representation theory of Clifford algebras, the space $\Gamma(S)$ of the smooth sections of $S$ is naturally bi-graded:
\[ \Gamma(S)= \Gamma(S)^+\oplus \Gamma(S)^-.\]
In the most interesting cases the Dirac operator is an odd operator: $\KD=\left(\begin{smallmatrix}0&\KD_-\\ \KD_+&0\end{smallmatrix}\right)$. It may be self-adjoint or skew self-adjoint, depending on the Clifford relation that defines the Clifford algebra. Because we will rely on the results of Alekseev and Meinrenken, we shall follow their convention, under which we have:
\begin{equation}
XY+YX = \langle X,Y\rangle 
\label{eq:clffrel}
\end{equation}
for vectors $X$ and $Y$ in the Euclidean space with inner product $\langle \ {,} \ \rangle $ that generate the Clifford algebra. This makes $\KD$  skew self-adjoint. Then the spectrum of $\KD$ is an unbounded discrete subset of $i\R$, each eigenvalue occurring with finite multiplicity. The space $\Gamma^2(S)$ of square-integrable sections of $S$ admits a Hilbert space direct sum decomposition into the eigenspaces of $\KD$. (For general reference on the analytic properties of Dirac operators, see Roe~\cite{roe}*{Chs.~5 \& 7}.) One can then build the heat diffusion operator $e^{t\KD^2}$ on $\Gamma^2(S)$ for $t\in(0,\infty)$. The main reason for taking this step arises from the \emph{McKean-Singer formula}, which relates the graded index of $\KD$ with the super trace of the heat diffusion operator; namely, 
\[ \Ind\KD = \Str(e^{t\KD^2}), \]
where $\Ind\KD$ is the difference between the dimensions of $\ker\KD_+$ and $\ker\KD_-$, and $\Str(e^{t\KD^2})$ is the difference between the trace of $e^{t\KD^2}$ on the even domain and that on the odd domain. The formula can be rewritten in terms of the \emph{heat kernel} associated with $\KD$, that is, the integral kernel $P_t$ of $e^{t\KD^2}$; it is a section of the external tensor product $S\boxtimes S^*$ over $M\times M$ such that 
\[
(e^{t\KD^2}\xi)(x) = \int_M P_t(x,y)\xi(y)\,\vol_y 
\]
for any $\xi$ in $\Gamma(S)$. Here $\vol_y$ is the value of the Riemannian volume form of $M$ at $y$. In terms of the heat kernel, the Mckean-Singer formula can be set down as:
\begin{equation}
\Ind\KD = \int_M \Str(P_t(x,x))\,\vol_x,
\label{eq:mckeansingerint}
\end{equation} 
where $\Str(P_t(x,x))$ calculates the super trace of $P_t(x,x)$ on the fiber $S_x$. For this reason, $x\mapsto\Str(P_t(x,x))\,\vol_x$ is referred to as the \emph{index density} on $M$. The restriction of the heat kernel onto the diagonal $(x,x)\in M\times M$ admits an asymptotic expansion for $t\rightarrow 0+$ of the following form (Seeley~\cite{seeley}): 
\[ P_t(x,x) \sim t^{-n/2}(a_0(x) + a_1(x)t+a_2(x)t^2+\dotsb), \]
where $n=\dim(M)/2$. Thus,
\[ \Str(P_t(x,x)) \sim t^{-n/2}(\Str(a_0(x)) + \Str(a_1(x))t+\Str(a_2(x))t^2+\dotsb).\]
The \emph{local index theorem} then states that the coefficients $\Str(a_{j}(x))$ with $0\le j< n/2$ are zero so that 
\[ \Str(P_t(x,x))\,\vol_x = \Str(a_{n/2}(x))\,\vol_{x}+O(t);\]
moreover, the leading term can be explicitly calculated in terms of the characteristic forms of $M$. The local index theorem yields the global index theorem via Equation~\eqref{eq:mckeansingerint}; this line of argument is often referred to as the ``heat kernel proof'' for the Atiyah-Singer index theorem.

Now suppose $M=G/K$. We examine the local index theorem for $G/K$ in the following setting. Endow a bi-invariant metric on $G$. Let $\mathfrak{g}$ and $\mathfrak{k}$ be the Lie algebras of $G$ and $K$, respectively,  and let $\mathfrak{p}$ be the orthogonal complement of $\mathfrak{k}$ in $\mathfrak{g}$. Let $E$ denote a finite-dimensional $K$-vector space that admits a $K$-equivariant module structure over the Clifford algebra $\Cl(\mathfrak{p})$ generated by $\mathfrak{p}$. We assume that $E$ is a graded $\Cl(\mathfrak{p})$-module with respect to the natural bi-grading on $\Cl(\mathfrak{p})$. We consider the homogeneous vector bundle $G\times_KE$. 

Our method for proving the local index theorem in the above setting is as follows.  The space $\Gamma(G\times_KE)$ is linearly isomorphic to the space $(C^{\infty}(G)\otimes E)^K$ of smooth $K$-invariant functions on $G$ with values in $E$. A $G$-equivariant differential operator on $(C^\infty(G)\otimes E)^K$ can be identified as a member of $(\mathcal{U}(\mathfrak{g})\otimes\Cl(\mathfrak{p}))^K$, where $\mathcal{U}(\mathfrak{g})$ is the universal enveloping algebra generated by $\mathfrak{g}$. This algebra, $(\mathcal{U}(\mathfrak{g})\otimes\Cl(\mathfrak{p}))^K$, is known as the \emph{relative Weil algebra} for the pair $(\mathfrak{g},\mathfrak{k})$; it was the main object of study in  \mbox{Alekseev} and \mbox{Meinrenken}'s work~\cite{alekmein2}. There they demonstrate that the algebraic structure of the relative Weil algebra singles out an element $\KD(\mathfrak{g},\mathfrak{k})$, which happens to be the cubic Dirac operator. We consider the Dirac operator $\KD_{\mathfrak{g}/\mathfrak{k}}$ on $\Gamma(G\times_KE)$ such that the following diagram commutes: 
\[
 \xymatrix@R=3ex@C=5em{ (C^\infty(G)\otimes E)^K\ar[r]^{\KD(\mathfrak{g},\mathfrak{k})} & (C^\infty(G)\otimes E)^K \\
	\Gamma(G\times_KE) \ar[u]^-{\cong}\ar[r]_{\KD_{\mathfrak{g}/\mathfrak{k}}} & \Gamma(G\times_KE) \ar[u]_-{\cong}
	} 
\]
This Dirac operator, $\KD_{\mathfrak{g}/\mathfrak{k}}$, is \emph{not} the Riemannian Dirac operator, which was used for instance in \mbox{Getzler's} proof~\cite{getzler} of the local index theorem. Our strategy is to deduce the heat kernel associated with $\KD_{\mathfrak{g}/\mathfrak{k}}$ from that of $\KD(\mathfrak{g},\mathfrak{k})$. 

To that end, we use the \emph{quantization map} 
\begin{equation*}
 \mathcal{Q}\colon (S(\mathfrak{g})\otimes\wedge(\mathfrak{p}))^K \to (\mathcal{U}(\mathfrak{g})\otimes\Cl(\mathfrak{p}))^K\label{eq:relweilsysy}
\end{equation*}
introduced by Alekseev and Meinrenken~\cite{alekmein}*{$\S$~6}; here $S(\mathfrak{g})$ is the symmetric algebra generated by $\mathfrak{g}$, and $\wedge(\mathfrak{p})$ is the exterior algebra generated by $\mathfrak{p}$. This map is a linear isomorphism that exerts a graded symmetrization with respect to certain generators. There is a natural identification of members of   $(S(\mathfrak{g})\otimes\wedge(\mathfrak{p}))^K$ with differential operators on $(C^\infty(\mathfrak{g})\otimes E)^K$. Let $\mathcal{L}$ be the preimage of $\KD(\mathfrak{g},\mathfrak{k})^2$ under the quantization map so that
\begin{equation*}
\mathcal{Q}(\mathcal{L})=\KD(\mathfrak{g},\mathfrak{k})^2.
\label{eq:kostantdiracqtmp}
\end{equation*}
The differential operator $\mathcal{L}$ is in fact a member of the subalgebra $S(\mathfrak{g})^G$. On this subalgebra the quantization map is identical to the Duflo isomorphism~\cites{duflo,alekmein}, which is an \emph{algebra} isomorphism. From this fact we deduce that the asymptotic expansion of the heat kernel $R_t$ of $\KD(\mathfrak{g},\mathfrak{k})^2$ must satisfy (Theorem~\ref{thm:asympkdhtk}):
\begin{equation}
 R_t^{\exp}= h_t^\mathfrak{g} \, \frac{1}{j}+O(t). \label{eq:asymptg}
\end{equation}
Here $R_t^{\exp}$ denotes the function $x\mapsto R_t(e,x)$ under the exponential chart, $h_t^\mathfrak{g}$ is the Gaussian function, and $j$ is the positive square root of the Jacobian determinant of the exponential map. What remains is to take the asymptotic equation~\eqref{eq:asymptg} and extract from it the asymptotic expansion of the heat kernel $P_t$ of $\KD^2_{\mathfrak{g}/\mathfrak{k}}$. For this task, we follow  the idea of \mbox{Berline} and \mbox{Vergne}~\cite{berlinevergne}, namely, that $P_t(xK,xK)$ is, as an element of $\End(E)$, equal to 
\[
 \int_K R_t(x,xh^{-1})\nu(h)\,dh, 
\]
where $\nu$ denotes the representation of $K$ on $E$, and $dh$ is the Haar measure on $K$. We shall see that throughout the proof almost all calculations are brought down to the level of Lie algebras owing to homogeneity.

In the next section, we describe the notations and conventions that we implement and the precise means by which we identify an element of the quantum Weil algebra as an equivariant differential operator on $(C^\infty(G)\otimes E)^K$. 

In the third section, using the algebraic properties of the cubic Dirac operator concatenated with the quantum Weil algebra, we derive the asymptotic expansion of the heat kernel of the (scalar) Laplacian on $(C^\infty(G)\otimes E)^K$; this is where the advantage of our method stands out.

And finally, in the last section, we deduce the local index theorem on $G/K$. The argument proceeds in the manner of Berline and Vergne~\cite{berlinevergne}, and we present only those calculations that are not mere adaptations of their work.

\section{Equivariant Differential Operators and the Relative Weil Algebra}\label{sec:equivdiffop}

Let $(G,K)$ be a pair of connected Lie groups where $K$ is a closed subgroup of $G$. Let $(\mathfrak{g},\mathfrak{k})$ be the associated Lie algebras. Since we are assuming that $G$ is connected, the existence of a bi-invariant metric on $G$ is equivalent to the existence of an inner product $\langle \ {,} \ \rangle$ on $\mathfrak{g}$ that is invariant under the adjoint representation $\ad\colon \mathfrak{g}\to \End(\mathfrak{g})$ so that 
\[ \langle \ad(X)Y,Z\rangle =- \langle Y,\ad(X)Z\rangle  \]
for any vectors $X$, $Y$,  and $Z$ in $\mathfrak{g}$. We shall assume such an inner product exists, as it is the case when $G$ is compact. 

Owing to the invariance of the inner product, the orthogonal decomposition $\mathfrak{g}=\mathfrak{k}\oplus\mathfrak{p}$ is preserved under the adjoint action of $\mathfrak{k}$. In particular, we have a Lie algebra homomorphism:
\[
\ad^{\mathfrak{p}}\colon \mathfrak{k} \to \so(\mathfrak{p}).
\]
Now the Lie algebra $\so(\mathfrak{p})$ is isomorphic to the Lie algebra $\spin(\mathfrak{p})$ of $\Spin(\mathfrak{p})$. So we have a Lie algebra homomorphism
\[\gamma^{\mathfrak{p}}\colon \mathfrak{k}\to\spin(\mathfrak{p}) \subset\Cl(\mathfrak{p}),
\]
which is characterized by the relation (this is due to Kostant \cite{kostantcliff}*{$\S$~2})
\[
 [\gamma^{\mathfrak{p}}(X),Y]=[X,Y]_{\mathfrak{g}}
 \]
for any $X\in\mathfrak{k}$ and any $Y\in \mathfrak{p}$;  the bracket on the left denotes the graded commutator in $\Cl(\mathfrak{p})$;  on the right is the Lie bracket in $\mathfrak{g}$. In terms of an orthonormal basis $\{ Y_i\}_{i=1}^{\dim\mathfrak{p}}$ for $\mathfrak{p}$, we have:
\begin{equation}
\gamma^{\mathfrak{p}}(X)=-\frac{1}{2}\sum_{i,j=1}^{\dim\mathfrak{p}}\langle X,[Y_i,Y_j]_{\mathfrak{g}}\rangle Y_iY_j. 
\label{eq:decmopsospin2}
\end{equation}

Throughout this article, we let $E$ denote an arbitrary $\Cl(\mathfrak{p})$-module that admits a $K$-action $\nu\colon K\to \Aut(E)$. We assume that the $\Cl(\mathfrak{p})$-multiplication on $E$ is $K$-equivariant, that is, the equation
\begin{equation}
 \nu(h)\circ Y \circ \nu(h)^{-1} = \Ad(h)(Y)
\label{eq:cliffkequiv}
\end{equation}
holds in $\End(E)$ for any $h\in K$ and any $Y\in\mathfrak{p}$. Here $\Ad$ denotes, as usual, the adjoint representation of $G$ on $\mathfrak{g}$. The induced Borel mixing space 
\[ E(G)\ldef G\times_KE \]
is the $K$-orbit space of $G\times E$ where $h\in K$ acts on $(g,v)\in G\times E$ by $(g,v)\cdot h=(gh,  \nu(h)^{-1} v)$. We denote the $K$-orbit of $(g,v)$ by $[g,v]$. The canonical projection $E(G)\to G/K$, $[g,v]\mapsto gK$, is a vector bundle with fiber $E$. A section of $E(G)$ can be identified with a function in $C^\infty(G)\otimes E$ that is $K$-invariant relative to $r\otimes \nu$, where $r$ denotes the right-regular action of $K$ on $C^\infty(G)$. (For details, see Kobayashi and Nomizu~\cite{kobayashinomizu1}*{Ex.~5.2, p.~76} or Berline, Getzler, and Vergne~\cite{bgv}*{Prop.~1.7, p.~19}.) This identification yields a linear isomorphism 
\begin{equation}
  \eta:\Gamma(E(G)) \xrightarrow{\sim}  (C^\infty(G)\otimes E)^K. \label{eq:equivfuncs}
\end{equation}

Our aim is to identify the elements of the relative Weil algebra 
\[ \qW(\mathfrak{g},\mathfrak{k}) \ldef  (\mathcal{U}(\mathfrak{g})\otimes\Cl(\mathfrak{p}))^K \]
with $G$-invariant differential operators on $(C^\infty(G)\otimes E)^K$. Since $E$ comes with a $\Cl(\mathfrak{p})$-module structure, we only need to specify the action of $\mathcal{U}(\mathfrak{g})$ on $C^\infty(G)$. For that purpose, it is sufficient to describe how the generating set $\mathfrak{g}$ acts on $C^\infty(G)$. Following \mbox{Alekseev} and \mbox{Meinrenken}, if $X\in\mathfrak{g}$ takes the role of a generator for $(\mathcal{U}(\mathfrak{g})\otimes\{1\})\cap\qW(\mathfrak{g},\mathfrak{k})$ then we shall denote it as:
\[ \hat X = X\otimes 1.\]
The natural step is to let $\hat{X}$ act on $C^\infty(G)$ as the left-invariant vector field on $G$ generated by $X$:
\begin{equation*} 
(\hat X\sigma)(g) = \left.\frac{d}{dt} {\sigma(g \exp(tX))} \right|_{0}
\label{eq:lfinvcfldact}
\end{equation*}
for $\sigma\in (C^\infty(G)\otimes E)^K$ and $g\in G$. Then we have:
\[ [\hat X,\hat Y]\sigma = \widehat{[X,Y]}_{\mathfrak{g}}\sigma.\]
In this way, we have an algebra isomorphism from $\mathcal{U}(\mathfrak{g})$ to the algebra of left-invariant differential operators on $G$. This completes the description of the $\mathcal{W}(\mathfrak{g},\mathfrak{k})$-action on $(C^\infty(G)\otimes E)^K$. 

Next, we demand that $\mathcal{W}(\mathfrak{g},\mathfrak{k})$ acts on $\Gamma(E(G))$ in a compatible way so that the following diagram is commutative (it is not a commutative diagram in the conventional sense):
\begin{equation*}\begin{split}
 \xymatrix@C=4em{ (C^\infty(G)\otimes E)^K\ar[r]^{\qW(\mathfrak{g},\mathfrak{k})} & (C^\infty(G)\otimes E)^K \\
	\Gamma(E(G)) \ar[u]^{\eta}_{\cong}\ar[r]_{\qW(\mathfrak{g},\mathfrak{k})} & \Gamma(E(G)) \ar[u]_{\eta}^{\cong}
	}
\end{split} \label{diag:diffophmlft}
\end{equation*}
Of particular interest is the resulting action of $\hat{Y}$ on
$\xi\in\Gamma (E(G))$ where $Y\in \mathfrak{p}$. 
It is well-known that, for any $g\in G$, the map
\begin{equation}
\label{eq:grplcldiffeohms}
\begin{array}{c@{\quad}c@{\quad}c}
\mathfrak{p} & \to &G/K,
\\
X&\mapsto &g\exp (X)K,
\end{array}
\end{equation}
is a local diffeomorphism near the origin (see
Helgason~\cite[Lem.~4.1, p.~123]{helgason}). It is then
straight forward to check that
\begin{equation*}
(\hat{Y}\xi)(gK) = \partial _{Y}\xi,
\end{equation*}
where $\partial _{Y}\xi$ denotes the directional derivative of $\xi$  
with respect to $Y$ under the chart provided by the local
diffeomorphism~\eqref{eq:grplcldiffeohms}; in other words, 
$(\hat{Y} \xi )(gK)= [g,\left.\frac{d}{dt}w(t)\right|_{t=0}]$,
where $t\mapsto w(t)$ is a curve in $E$ defined on some interval containing $0$ such that
$\xi (g\exp (tY)K) = [g\exp (tY),w(t)]$.
More generally, we have
\begin{equation*}
(\hat{Y}\otimes v\xi)(gK)  = v\partial _{Y}\xi
\end{equation*}
for any $v\in \operatorname{Cl}(\mathfrak{p})$.

\begin{rmk}\label{rmk:equivweilact}
The derivative $\partial$ is, in fact, a covariant derivative. This type of covariant derivative appears in the literature  (not always in the same setting as ours) with various premodifiers. Following Kobayashi and Nomizu  \cite{kobayashinomizu2}, we shall refer to $\partial$ as the \emph{canonical} covariant derivative. (This covariant derivative is also \emph{reductive} in the sense that its torsion and curvature are parallel; see  \cite{kobayashinomizu2}*{Thm.~2.6, p.~193} and \cite{kostantred}*{Thm.~1, p.~37}.) 

The connection $1$-form $\theta$ associated with $\partial$ is as follows: Let $\omega$ be the (left-invariant) Maurer-Cartan form on $G$ so that for a left-invariant vector field $\tilde X$ on $G$ generated by $X\in\mathfrak{g}$ we have $\omega(\tilde X)=X$. Let $\pr_{\mathfrak{k}}$ denote the orthogonal projection of $\mathfrak{g}$ onto $\mathfrak{k}$. Then $\theta$ is the $\mathfrak{k}$-valued $1$-form ${\pr_{\mathfrak{k}}}\circ\omega$. 

Note that the canonical connection $\theta$ characterizes tangent vectors on $G$ that are orthogonal to $K$-orbits as horizontal vectors; in other words, the left-translations of $\mathfrak{p}$ are precisely the horizontal subspaces for the tangent spaces of $G$, and thus, $\mathfrak{p}$ is identified with the tangent space of $G/K$ at any position by means of the local diffeomorphism~\eqref{eq:grplcldiffeohms}.
\end{rmk}

\section{The Asymptotic Heat Kernel on Compact Lie Groups}

The algebraic structure of $\qW(\mathfrak{g},\mathfrak{k})$ picks out a distinguished element $\KD(\mathfrak{g},\mathfrak{k})$, namely, the cubic Dirac operator. To wit, $\qW(\mathfrak{g},\mathfrak{k})$ has a differential graded algebra structure of which the differential is inner with respect to $\mathcal{D}(\mathfrak{g},\mathfrak{k})$ (Alekseev and Meinrenken~\cite{alekmein2}*{Prop.~6.4, p.~321}). The precise form of $\KD(\mathfrak{g},\mathfrak{k})$ is:
\begin{equation*} \KD(\mathfrak{g},\mathfrak{k}) = \sum_{i=1}^{\dim\mathfrak{p}} Y_i\otimes Y_i +1\otimes \frac{1}{3} \sum_{i=1}^{\dim\mathfrak{p}}Y_i\gamma^\mathfrak{p}(Y_i) \quad \in \qW(\mathfrak{g},\mathfrak{k}) = (\mathcal{U}(\mathfrak{g})\otimes\Cl(\mathfrak{p}))^K, \label{eq:kostantdirachmsp}
\end{equation*}
where $\{Y_i\}_{i=1}^{\dim\mathfrak{p}}$ is any orthonormal basis for $\mathfrak{p}$, and $\gamma^{\mathfrak{p}}(Y_i)$ is defined using the same formula expressed in Equation~\eqref{eq:decmopsospin2}.

The objective of this section is to derive the asymptotic heat kernel of $\KD(\mathfrak{g},\mathfrak{k})^2$ on $(C^\infty(G)\otimes E)^K$ utilizing its algebraic properties associated with $\qW(\mathfrak{g},\mathfrak{k})$, especially those related to the quantization map
\[ \mathcal{Q}\colon (S(\mathfrak{g})\otimes\wedge(\mathfrak{p}))^K \to (\mathcal{U}(\mathfrak{g})\otimes\Cl(\mathfrak{p}))^K.
\]
We need not know all the details of this map except that $\mathcal{Q}$, upon restriction, yields the Duflo isomorphism $S(\mathfrak{g})^G\cong \mathcal{U}(\mathfrak{g})^G$, which is an \emph{algebra} isomorphism~\cite{alekmein}*{Thm.~7.1, p.164}. 

Of course, we need to relate the above algebraic fact with differential operators. To that end, we can define the action of $(S(\mathfrak{g})\otimes\wedge(\mathfrak{p}))^K$ on $E$-valued functions on $\mathfrak{g}$ in a similar manner as done in Section~\ref{sec:equivdiffop}, which identifies $S(\mathfrak{g})^G$ with the algebra of $G$-invariant constant coefficient differential operators on $\mathfrak{g}$. Then $\mathcal{Q}$ relates an operator in $S(\mathfrak{g})^G$ with a bi-invariant operator on $G$ by comparing the full symbols at the origin in $\mathfrak{g}$ and the identity in $G$, respectively. In this picture, the standard Laplacian $\Delta_{\mathfrak{g}}$ on $\mathfrak{g}$ is a member of $S(\mathfrak{g})^G$. We shall see that $\mathcal{Q}(\Delta_{\mathfrak{g}})/2$ (the factor $1/2$ is due to our Clifford relation~\eqref{eq:clffrel}) is a differential operator whose principal symbol is the same as $ \KD(\mathfrak{g},\mathfrak{k})^2$.  In fact, if $E=S\otimes V$ where $S$ is a graded irreducible $\Cl(\mathfrak{p})$-module and $V$ is an irreducible $\mathfrak{k}$-vector space, then (Theorem~\ref{thm:genlapfml}):
\[ \KD(\mathfrak{g},\mathfrak{k})^2=\frac{1}{2} \mathcal{Q}(\Delta_{\mathfrak{g}}) + \text{constant}, \]
where the constant depends on the $\mathfrak{k}$-action on $V$. 

\subsection{Some Terminology and Conventions}\label{sec:htknl}
Just like any other Dirac operator the square of $\KD(\mathfrak{g},\mathfrak{k})$ is a generalized Laplacian. By a generalized Laplacian we mean the following. Let $F$ be a vector bundle over a compact oriented Riemannian manifold $M$. We assume that $F$ is equipped with a fiber-wise inner product in a smooth manner. A generalized Laplacian on $\Gamma(F)$ is a differential operator $L$ whose expression in local coordinates takes the form of:
\begin{equation}
 L=\frac12\sum_{i,j=1}^{\dim M}\beta^{ij}\partial_i\partial_j + (\mathrm{lower \ order \ part}), \label{eq:genlapdef}
\end{equation}
where $\beta^{ij}$ is the $(i,j)$-entry of the inverse of the matrix $[\beta_{ij}]$ coming from the metric $\beta = \sum_{i,j} \beta_{ij}dx^idx^j$ on $M$ (we are distinguishing the upper indices from the lower indices). The $1/2$ factor in Equation~\eqref{eq:genlapdef} is an insignificant part in our discussion; it is placed here on account of the Clifford relation~\eqref{eq:clffrel}, which is the convention used by Alekseev and Meinrenken~\cite{alekmein}. The heat kernel $P_t$ of $L$, parametrized by ``time'' $t$, is a smooth section of the bundle $F\boxtimes F^*$ such that
\[ (e^{tL}\xi)(x) = \int_M P_t(x,y)\xi(y)\,\vol_y \]
for any $t>0$ and any $\xi\in \Gamma(F)$. For the existence of  $P_t$ and its fundamental properties, see Roe~\cite{roe}*{Prop.~5.31, p.~83; Prop.~7.5, p.~96}.

Suppose $M$ is a compact homogeneous space $G/K$. We adopt the following notation for the left coset of $K$ that contains $x\in G$: 
\[\bar x\ldef  xK.\]
Assume that the bundle $F$ is homogeneous over $G/K$, in the sense that 
\[ F = V(G)=G\times_KV \]
for some finite-dimensional $K$-vector space $V$. There is a natural $G$-action on $F$ defined by
\[ \alpha_{g}\colon [g',v] \mapsto [gg',v]  \]
for $g\in G$ and $[g',v]\in V(G)$. This induces a $G$-action on the sections of $F$, namely, for $g\in G$, $\xi\in \Gamma(F)$, and $\bar x\in G/K$, 
\begin{equation}
 (g\cdot \xi)(\bar x) \ldef  \alpha_g(\xi(g^{-1} \bar x)),
\label{eq:gactsec}
\end{equation}
where $g^{-1}\bar x\ldef  g^{-1} xK$. 

Suppose the generalized Laplacian $L$ on $\Gamma(F)$ is equivariant with respect to the action~\eqref{eq:gactsec}. Then the heat diffusion operator $e^{tL}$ is also equivariant. Provided that the metric on $M$ is $G$-invariant, the heat kernel $P_t$ is completely determined by the function
\[ 
p_t\colon \bar y \mapsto  P_t(\bar e,\bar y). 
\]
We shall refer to $p_t$ as the \emph{heat convolution kernel} of $L$.  We may view it as a function $
(0,\infty)\to \Gamma(F_{\bar e}\boxtimes F^*)$, $t\mapsto p_t$. 

For small time $t$, the heat convolution kernel behaves like  the Gaussian kernel\footnote{This is the heat kernel for $1/2$ times the Laplacian; the heat kernel for the usual Laplacian can be obtained by the substitution $t\mapsto 2t$ in \eqref{eq:rscgaussiankernel}.}
\begin{equation}
  h_t(\bar x)=\frac{e^{-d(\bar x)^2/2t}}{(2\pi t)^{\dim M/2}}, \label{eq:rscgaussiankernel}
\end{equation}
where $d(\bar x)$ denotes the distance between $\bar e$ and $\bar x$; stating more precisely, $p_t$ admits an asymptotic expansion of the form 
\begin{equation}
 p_t \sim h_t (a_0+a_1t+a_2t^2+\dotsb) \label{eq:asympexpg}
\end{equation}
for $t\to0+$, valid in the Banach space of  $C^r$-sections of $F_{\bar e}\boxtimes F^*$ for all non-negative integer $r$ (see Roe~\cite{roe}*{Thm.~7.15, p.~101}). The right-hand side of~\eqref{eq:asympexpg}, often referred to as the \emph{asymptotic heat kernel} for $p_t$, is obtained as the formal solution to the heat equation associated with $L$, namely, 
\[
 (\partial_t+L)p_t = 0.\label{eq:diffeqasympexp}
\]
That is, we formally solve the heat equation by substituting $p_t$ with $h_t\sum^\infty_{i=0}a_it^i$, under the condition that $a_0(\bar e)$ is the identity operator on the fiber $F_{\bar e}$. For details, see Roe~\cite{roe}*{Thm.~7.15, p.~101} or Berline et al.~\cite{bgv}*{Thm.~2.26, p.~83; Thm.~2.29, p.~85}.

\subsection{A Few Words on the Duflo Isomorphism}

As mentioned in the beginning of this section, it is decisive in our calculation that the restriction of the  quantization map $\mathcal{Q}$ to $S(\mathfrak{g})^G$ is the Duflo isomorphism. It is also necessary to translate the algebraic definition of $\mathcal{Q}$ into the language of differential operators and distributions. 

To that end, let $\mathcal{E}'(\mathfrak{g})^G$ and $\mathcal{E}'(G)^G$ denote the algebra of compactly supported $G$-invariant distributions on $\mathfrak{g}$ and $G$, respectively. Let $j_\mathfrak{g}$ be the analytic function on $\mathfrak{g}$ defined by: 
\[ j_\mathfrak{g}(X) = \det\nolimits^{1/2}\biggl[\frac{\sinh(\ad(X/2))}{{\ad}(X/2)}\biggr]. \]
As it is well-known, $j_\mathfrak{g}^2(X)$ calculates the Jacobian determinant of the exponential map at $X$ if the exponential map is diffeomorphic near $X$. Now consider the map 
\begin{equation}
 {\exp_*}\circ j\colon \mathcal{E}'(\mathfrak{g})^G\to \mathcal{E}'(G)^G,
\label{eq:dufdist}
\end{equation}
where $\exp_*$ denotes the push-forward along the exponential map, and $j$ denotes the multiplication by $j_\mathfrak{g}$. Though it is an abuse of notation, we shall denote the map~\eqref{eq:dufdist} also by $ \mathcal{Q}$, for $\exp_{*}\circ j$ does agree with the quantization map when invariant differential operators are identified as point-supported distributions \cite[Thm.~6.1, p.~159]{alekmein}. Then, for any differential operator $P \in S(\mathfrak{g})^G$, its image $\mathcal{Q}(P)\in \mathcal{U}(\mathfrak{g})^G$ under the quantization map  is the bi-invariant differential operator on $G$ that satisfies
\begin{equation}
 \mathcal{Q}(P)\mathcal{Q}(\delta_0) = \mathcal{Q}(P\delta_0),
\label{eq:dufdistdef}
\end{equation}
where $\delta_0$ is the Delta distribution on $\mathfrak{g}$ supported at the origin. (Accordingly $\mathcal{Q}(\delta_0)$ is the Delta distribution $\delta_e$ on $G$ supported at the identity.) 

It is known that the map~\eqref{eq:dufdist} is an algebra homomorphism, that is,
\begin{equation} 
\mathcal{Q}(u)*\mathcal{Q}(v)=\mathcal{Q}(u*v) \label{eq:kveq}
\end{equation}
for any $u$ and $v$ in $\mathcal{E}'(\mathfrak{g})^G$. This is the substance of the Kashiwara-Vergne conjecture \cite{kv}, for which there are now several proofs for any Lie group $G$. In our case, where $G$ is compact, it is easy to verify Equation~\eqref{eq:kveq} using the \mbox{Kirillov} character formula (Dooley and Wildberger~\cite{dooleywild}). But there are proofs that do not involve representation theory or the structure theory of the group $G$. In particular, when $\mathfrak{g}$ admits a nondegenerate symmetric bilinear form (like in our case), then Equation~\eqref{eq:kveq} is a ``trivial'' consequence (Alekseev and Torossian~\cite{alektor}) of the fact that the Campbell-Baker-Hausdorff formula can be written in the following form:
\[ \log(\exp(X)\exp(Y)) = X+Y -(e^{\ad(X)}-1)A(X,Y)-(1-e^{-\ad(Y)})B(X,Y).\]

\subsection{Calculating the Asymptotic Heat Kernel}
According to the representation theory of Clifford algebras, there is a unique (up to isomorphism) irreducible graded $\Cl(\mathfrak{p})$-module when $\mathfrak{p}$ is odd-dimensional. In the even-dimensional case there are two irreducible graded modules; as vector spaces they are both the direct sum of the positive and negative half-spinors, and the ambiguity can be removed by declaring the positive and negative subspaces as the even and odd subspaces, respectively (see \cite{encg}*{pp.~188--189}). Henceforth we may, regardless of the parity of the dimension of $\mathfrak{p}$, speak of \emph{the} graded irreducible spinor space $S$.

Having a unique irreducible graded $\Cl(\mathfrak{p})$-module, the equivariance condition~\eqref{eq:cliffkequiv} implies that there is a Lie algebra representation $\tau\colon \mathfrak{k} \to \End(V)$ such that we have the equality
\[ E= S\otimes V \]
as $\mathfrak{k}$-spaces, where $\mathfrak{k}$ acts on the right via $\gamma^{\mathfrak{p}}\otimes\Id+\Id\otimes\tau$ (see Goette~\cite{goette}*{Lem.~3.4, P.~25}). We shall simply write this fact as:
\begin{equation}
 \nu_* = \gamma^{\mathfrak{p}} + \tau. \label{eq:indlierepclifm}
\end{equation}
The space $V$ serves as an auxiliary space on which $\Cl(\mathfrak{p})$ acts trivially. Each Lie algebra representation  $\nu_*$, $\gamma^{\mathfrak{p}}$, and $\tau$ extends to $\mathcal{U}(\mathfrak{k})$ as an algebra homomorphism.

\emph{We shall assume that $V$ is $\mathfrak{k}$-irreducible with highest weight $\mu$.} (For reference on the theory of highest weight, we refer to Duistermaat and Kolk~\cite{duistermaatkolk}*{$\S$~4.9}.)  

\begin{lem}\label{lem:diagcaskinvsec}
Let $\diag_\qW$ denote the (diagonal) embedding $\mathcal{U}(\mathfrak{k})\hookrightarrow\qW(\mathfrak{g},\mathfrak{k})$ induced by $X\mapsto \hat X + \gamma^{\mathfrak{p}}(X)$. Let $\Omega_{\mathfrak{k}}$ denote the Casimir element in $\mathcal{U}(\mathfrak{k})$ so that $\Omega_{\mathfrak{k}}=\sum_{i=1}^{\dim\mathfrak{k}} X_iX_i$ where $\{X_i\}_{i=1}^{\dim\mathfrak{k}}$ is any orthonormal basis for $\mathfrak{k}$. As differential operators on $(C^\infty(G)\otimes E)^K$, we have 
\[ \diag_\qW \Omega_{\mathfrak{k}} = -\|\mu+\rho_{\mathfrak{k}}\|^2+\|\rho_{\mathfrak{k}}\|^2, \]
where $\rho_{\mathfrak{k}}$ is the Weyl vector of $\mathfrak{k}$, and $\|\cdot\|$ is the norm induced by the inner product on $\mathfrak{g}$.
\end{lem}
\begin{proof}
Let $\sigma \in (C^\infty(G)\otimes E)^K$. By the $K$-invariance of $\sigma$, we have
\[ \sigma(g\exp(X)) = e^{-\nu_*(X)}\sigma(g) \]
for every $g\in G$ and every $X\in\mathfrak{k}$. Differentiating both sides with respect to $X$, we get
\[ \hat X\sigma = -\nu_*(X)\sigma.\]
By Equation~\eqref{eq:indlierepclifm}, the right-hand side is $-(\gamma^{\mathfrak{p}}(X)+\tau(X))\sigma$. And, by definition, $\diag_\qW(X)=\hat X+\gamma^{\mathfrak{p}} (X)$. So the above equation can be rewritten as
\[ \diag_\qW(X)\sigma = -\tau(X)\sigma.\]
Hence, $\diag_\qW=-\tau$; this continues to hold as an equality between  algebra homomorphisms on $\mathcal{U}(\mathfrak{k})$;  so
\[ \diag_\qW(\Omega_{\mathfrak{k}}) =  \tau(\Omega_{\mathfrak{k}}) \]
as differential operators on  $(C^\infty(G)\otimes E)^K$. The right-hand side is a constant operator, owing to the irreducibility of $\tau$ and Schur's lemma. Its value is well-known to be $-{\|\mu+\rho_{\mathfrak{k}}\|^2}+\|\rho_{\mathfrak{k}}\|^2$ (see Kostant~\cite{kostant}*{Rmk.~1.89, p.~469}).
\end{proof}

\begin{thm}\label{thm:genlapfml}
As differential operators on $(C^\infty(G)\otimes E)^K$, we have
\begin{equation*}
\KD(\mathfrak{g},\mathfrak{k})^2 =   \frac12\Omega_{\mathfrak{g}} +\frac12\bigl( \|\mu+\rho_{\mathfrak{k}}\|^2 -  \|\rho_{\mathfrak{g}}\|^2\bigr).
	  \label{eq:genlapfml}
\end{equation*}
\end{thm}
\begin{rmk}
A more general formula involving a $1$-parameter family of Dirac operators can be found in Goette~\cite{goette}*{Lem.~1.17, p.~9}; see also Agricola~\cite{agricola}*{Thm.~3.3, p.~550; Lem.~3.6, p.~552}. In both cases there is a difference by a factor of $-1/2$ from our formula because they use the Clifford relation $XY+YX=-2\langle X,Y\rangle $ instead of \eqref{eq:clffrel}.
\end{rmk}
\begin{proof}
In Kostant~\cite{kostant}*{Prop.~1.84, p.~468; Thm.~2.13, p.~474}, the following formulas are derived: 
\begin{equation}
 \KD(\mathfrak{g},\mathfrak{k})^2 =\frac{1}{2}\Omega_{\mathfrak{g}} -  \frac{1}{2} \diag_{\qW} \Omega_{\mathfrak{k}} +\frac{1}{48}\tr_{\mathfrak{g}}\Omega_{\mathfrak{g}} -\frac{1}{48}\tr_{\mathfrak{k}}\Omega_{\mathfrak{k}},
 \label{eq:sqrreldirac}
\end{equation}
and
\begin{equation}
\frac{1}{24}\tr_{\mathfrak{g}}\Omega_{\mathfrak{g}} = -\|\rho_{\mathfrak{g}}\|^2.
\label{eq:kostantcasimirtrfml}
\end{equation}
Combining the two equations, we have
\[ \KD(\mathfrak{g},\mathfrak{k})^2 =\frac12\Omega_{\mathfrak{g}} -  \frac12 \diag_{\qW} \Omega_{\mathfrak{k}} -\frac{1}{2}\|\rho_{\mathfrak{g}}\|^2+\frac12\|\rho_{\mathfrak{k}}\|^2. 
\]	 
The assertion now follows from Lemma~\ref{lem:diagcaskinvsec}.
\end{proof}

\begin{thm}\label{thm:duffltlap}
Let $\Delta_{\mathfrak{g}}$ be the element in $S(\mathfrak{g})$ defined by $\Delta_{\mathfrak{g}}=\sum_{i=1}^{\dim\mathfrak{g}}\hat X_i\hat X_i$ where $\{X_i\}_{i=1}^{\dim\mathfrak{g}}$ is any orthonormal basis for $\mathfrak{g}
$ (the definition is independent of the choice of the basis). We have
\[ \mathcal{Q}(\Delta_{\mathfrak{g}})= \Omega_{\mathfrak{g}}-\|\rho_{\mathfrak{g}}\|^2.\]
\end{thm}
\begin{proof}
The relative Weil algebra $\qW(\mathfrak{g},\mathfrak{k})$ is a subalgebra of $\qW(\mathfrak{g})\ldef \qW(\mathfrak{g},\{0\})$. The quantization map for the pair $(\mathfrak{g},\mathfrak{k})$ is by definition the restriction of the quantization map for the pair $(\mathfrak{g},\{0\})$. According to Alekseev and Meinrenken~\cite{alekmein2}*{proof of Prop.~5.4, p.~319}, the element $\KD(\mathfrak{g})\ldef \KD(\mathfrak{g},\{0\})$ satisfies
\[ \KD(\mathfrak{g})^2= \frac{1}{2}\mathcal{Q}(\Delta_{\mathfrak{g}}).\]
Meanwhile, owing to Equation~\eqref{eq:sqrreldirac} and \eqref{eq:kostantcasimirtrfml} for the pair $(\mathfrak{g},\{0\})$, we have
\[ \KD(\mathfrak{g})^2= \frac{1}{2}\Omega_{\mathfrak{g}} - \frac{1}{2} \|\rho_{\mathfrak{g}}\|^2.\]
Hence, the equation
\[ \mathcal{Q}(\Delta_{\mathfrak{g}})=\Omega_{\mathfrak{g}}-\|\rho_{\mathfrak{g}}\|^2 \]
holds in $\qW(\mathfrak{g})$. But $\Delta_{\mathfrak{g}}$ is in $S(\mathfrak{g})^G$, so its image $\mathcal{Q}(\Delta_{\mathfrak{g}})$ is in $\mathcal{U}(\mathfrak{g})^G$. The Casimir element $\Omega_{\mathfrak{g}}$ is also in $\mathcal{U}(\mathfrak{g})^G$. Therefore, the above equation is valid in $\qW(\mathfrak{g},\mathfrak{k})$.
\end{proof}

\begin{rmk}
As differential operators, $\Delta_{\mathfrak{g}}$ and $\Omega_{\mathfrak{g}}$ are, respectively, the Laplacians on $\mathfrak{g}$ and $G$.
\end{rmk}

\begin{prop}\label{prop:duflapcptg2}
Let $\sigma$ be a smooth class function on $G$, and denote its pullback along the exponential map by $\sigma^{\exp}$. Then the equation
\[ (\mathcal{Q}(\Delta_{\mathfrak{g}})\sigma)^{\exp}=\Bigl(\frac{1}{j_\mathfrak{g}}\circ\Delta_{\mathfrak{g}}\circ j_\mathfrak{g} \Bigr)\sigma^{\exp} \]
holds on a sufficiently small neighborhood of the origin of $\mathfrak{g}$.
\end{prop}
\begin{rmk}
This is a special case of Theorem~1.12 of Helgason \cite{helgasoncbms}*{p.~30}. Still we present our proof here as it demonstrates that we need not any structure theory of $G$ for this result.
\end{rmk}
\begin{proof}
Let $f$ be a smooth function on $\mathfrak{g}$ that is $G$-invariant. Let $V$ be a neighborhood of the origin of $\mathfrak{g}$, sufficiently small so that the exponential map maps $V$ diffeomorphically onto its image. Using a $G$-invariant bump function if necessary, we may assume that $f$ is compactly supported in $V$. Then we can define $\Duf(f)\in C^\infty(G)$ that is supported in $\exp(V)$ such that the equation $\Duf(f)^{\exp} = f/j_\mathfrak{g}$ holds. This way we have
\begin{equation}
 \langle \Duf u, \Duf f \rangle = \langle u , f \rangle,
\label{eq:funcduf}
\end{equation}
for any $u\in\mathcal{E}'(\mathfrak{g})^G$. Under this notation, it is sufficient to prove that
\begin{equation}
\mathcal{Q}(\Delta_{\mathfrak{g}}) \Duf( f) = \Duf( \Delta_{\mathfrak{g}} f).
\label{eq:kvlap}
\end{equation} 
Indeed, for any class function $\sigma$, there is some $f$ such that $\sigma=\Duf(f)$ holds near the identity, and Equation~\eqref{eq:kvlap} implies that
\[ (\mathcal{Q}(\Delta_{\mathfrak{g}})\sigma)^{\exp} =(\mathcal{Q}(\Delta_\mathfrak{g})\mathcal{Q}(f))^{\exp}=\mathcal{Q}(\Delta_\mathfrak{g}f)^{\exp}= \frac{1}{j_\mathfrak{g}}\Delta_{\mathfrak{g}} f=\frac{1}{j_\mathfrak{g}}\Delta_{\mathfrak{g}}(j_\mathfrak{g}\sigma^{\exp}),\]
which is what we want.

To prove Equation~\eqref{eq:kvlap}, we only need the following fact: For any distributions $U$ and $V$ in $\mathcal{E}'(G)^G$, we have
\begin{equation}
 \mathcal{Q}(\Delta_{\mathfrak{g}})(U*V) = U*\mathcal{Q}(\Delta_{\mathfrak{g}})V.
\label{eq:dufdiff}
\end{equation}
This equation can be checked by simple calculations. Equations~\eqref{eq:dufdistdef}, \eqref{eq:kveq}, \eqref{eq:funcduf}, and \eqref{eq:dufdiff} can now be used to verify Equation~\eqref{eq:kvlap} as follows: For any test function $\phi$ on $G$,
\begin{align*}
&\langle \mathcal{Q}(\Delta_{\mathfrak{g}})\Duf(f), \phi \rangle
=\langle \mathcal{Q}(\Delta_{\mathfrak{g}})(\Duf(f)*\delta_e), \phi \rangle = 
\langle \Duf(f) * \mathcal{Q}(\Delta_{\mathfrak{g}})\delta_e, \phi \rangle \\
	&=\langle \mathcal{Q}(f)*\mathcal{Q}(\Delta_\mathfrak{g})\mathcal{Q}(\delta_0),\phi\rangle
	=  \langle \Duf(f) * \Duf(\Delta_{\mathfrak{g}}\delta_0), \phi \rangle 
	= \langle \Duf(f* \Delta_{\mathfrak{g}}\delta_0), \phi\rangle \\
	&= \langle \Duf(f* \Delta_{\mathfrak{g}}\delta_0), \Duf(j_\mathfrak{g}\phi^{\exp}) \rangle
	= \langle f * \Delta_{\mathfrak{g}}\delta_0, j_\mathfrak{g}\phi^{\exp} \rangle
	= \langle f\otimes \Delta_{\mathfrak{g}} \delta_0,j_\mathfrak{g}\phi^{\exp} \rangle \\
	&= \langle f,\Delta_{\mathfrak{g}}(j_\mathfrak{g}\phi^{\exp}) \rangle
	= \langle \Delta_{\mathfrak{g}}f,j_\mathfrak{g}\phi^{\exp} \rangle
	=\langle \Duf(\Delta_{\mathfrak{g}}f),\Duf(j_\mathfrak{g}\phi^{\exp}) \rangle\\
&=\langle \Duf(\Delta_{\mathfrak{g}}f),\phi \rangle.
\end{align*}
This proves Equation~\eqref{eq:kvlap}, and we are done.
\end{proof}

Proposition~\ref{prop:duflapcptg2} immediately leads to the asymptotic heat kernel of $\mathcal{Q}(\Delta_\mathfrak{g})/2$:
\begin{prop}\label{prop:asymhtkduflap2}
Let $q_t$ be the heat convolution kernel of $\mathcal{Q}(\Delta_{\mathfrak{g}})/2$. We have
\begin{equation}
 q_t^{\exp}\sim \frac{h_t^\mathfrak{g}}{j_\mathfrak{g}} 
 \label{eq:duflapahtk}
\end{equation}
for $t\to0+$, valid in a neighborhood of the origin. In other words, the coefficients in the asymptotic expansion~\eqref{eq:asympexpg} for this case are $a_0(x)=1/j^{\log}_\mathfrak{g}(x)$ and $a_n(x)=0$ for $n\ge1$.
\end{prop}

\begin{proof}
Let $s_t\ldef  h_t^\mathfrak{g}\sum_{i=0}^{\infty}a_it^i$ be the asymptotic heat kernel for $q_t^{\exp}$ (or $q_t^{\exp}$ multiplied by a bump function that has constant value $1$ on some neighborhood of the origin). The coefficients $a_i$ are obtained by acknowledging $s_t$ as the formal solution to the heat equation 
\[ \Bigl(\partial_t+\frac12\mathcal{Q}(\Delta_{\mathfrak{g}})\Bigr)q_t =0 \]
under the exponential chart. Because the convolution kernel $q_t$ is a class function (this follows from the fact that $\mathcal{Q}(\Delta_{\mathfrak{g}})$ is bi-invariant), the above differential equation under the exponential chart takes the form of (by Proposition~\ref{prop:duflapcptg2}): 
\[ \frac{1}{j_\mathfrak{g}}\Bigl(\partial_t-{\frac12}\Delta_{\mathfrak{g}}\Bigr)j_\mathfrak{g}q_t^{\exp}=0.\]
But the fundamental solution for the Euclidean heat operator $(\partial_t-\frac12\Delta_{\mathfrak{g}})$ is  the Gaussian $h_t^\mathfrak{g}$; hence $h_t^\mathfrak{g}/j_\mathfrak{g}$ is the formal solution for $q_t^{\exp}$ near the origin. This proves that $s_t=h_t^\mathfrak{g}/j_\mathfrak{g}$ is the asymptotic heat kernel for $q_t^{\exp}$.
\end{proof}

\begin{rmk}
\begin{enumerate}
\item The vanishing of the heat coefficients except for the lowest one is astonishing. So it would be worthwhile to give the following alternative proof. Let $\psi$ be a $G$-invariant bump function on $\mathfrak{g}$ centered at the origin. The function $\psi h_t^\mathfrak{g}$ defines a distribution in $\mathcal{E}'(\mathfrak{g})^G$ and, hence, is subject to $\mathcal{Q}$; the resultant in $\mathcal{E}'(G)^G$ is the distribution defined by the function $\mathcal{Q}( h_t^\mathfrak{g}\psi )$. Denote the integral operator on $(C^\infty(G)\otimes E)^K$ with kernel $\mathcal{Q}(h_t^\mathfrak{g}\psi)$ by $\mathcal{Q}(e^{t\Delta_{\mathfrak{g}}})$.  Then Proposition~\ref{prop:asymhtkduflap2} is equivalent to the asymptotic equality
\[  (e^{t\mathcal{Q}(\Delta_{\mathfrak{g}})/2})\sigma(e) \sim \mathcal{Q}(e^{t\Delta_{\mathfrak{g}}/2})\sigma(e)\]
for $t\to0+$. This asymptotic equality can be deduced using the fact that $\mathcal{Q}(\Delta_{\mathfrak{g}}^n)=\mathcal{Q}(\Delta_{\mathfrak{g}})^n$ as follows:
\begin{align*}
 (e^{t\mathcal{Q}(\Delta_{\mathfrak{g}})/2})\sigma(e) 
	&\sim \sum_{n=0}^\infty\Bigl< \delta_e, \mathcal{Q}(\Delta_{\mathfrak{g}})^n\sigma\Bigr>\frac{t^n}{n!2^n}
	= \sum_{n=0}^\infty\Bigl< \mathcal{Q}(\delta_0), \mathcal{Q}(\Delta_{\mathfrak{g}}^n)\sigma\Bigr>\frac{t^n}{n!2^n}\\
	&= \sum_{n=0}^\infty\Bigl< \mathcal{Q}(\Delta_{\mathfrak{g}}^n)\mathcal{Q}(\delta_0), \sigma\Bigr>\frac{t^n}{n!2^n}
	= \sum_{n=0}^\infty\Bigl< \mathcal{Q}(\Delta_{\mathfrak{g}}^n\delta_0), \mathcal{Q}(j_\mathfrak{g} \sigma^{\exp})\Bigr>\frac{t^n}{n!2^n}\\
	&= \sum_{n=0}^\infty\Bigl< \Delta_{\mathfrak{g}}^n\delta_0, j_\mathfrak{g} \sigma^{\exp}\Bigr>\frac{t^n}{n!2^n}
	=\sum^\infty_{n=0}\Bigl\langle\delta_0,\Delta_\mathfrak{g}^n(j_\mathfrak{g}\sigma^{\exp})\Bigr\rangle\frac{t^n}{n!2^n}\\
	&\sim  e^{t\Delta_{\mathfrak{g}}/2} (j_\mathfrak{g}\sigma^{\exp})(0)
	=\langle h_t^\mathfrak{g}, j_\mathfrak{g}\sigma^{\exp}\rangle \sim \langle   h_t^\mathfrak{g} \psi, j_\mathfrak{g}\sigma^{\exp}\rangle =\langle \mathcal{Q}(h_t^\mathfrak{g}\psi),\sigma \rangle\\
	&=\mathcal{Q}(e^{t\Delta_{\mathfrak{g}}/2})\sigma(e).
\end{align*}

\item From the asymptotic heat kernel~\eqref{eq:duflapahtk} we can deduce  the asymptotic expansion for the heat trace $Z(t)\ldef \tr(e^{t\Omega_\mathfrak{g}})$  of $G$. (Here the trace is over the square-integrable functions on $G$.) If $P_t$ is the heat kernel of $\Omega_\mathfrak{g}$, then $Z(t)=\int_GP_t(x,x)\,\vol_x=P_t(e,e)\,\vol(G)$ (the last equality owes to the $G$-invariance of $\Omega_\mathfrak{g}$). Now invoking Theorem~\ref{thm:duffltlap} and Proposition~\ref{prop:duflapcptg2}, we conclude that
\[ Z(t) \sim \frac{\vol(G)}{(4\pi t)^{\dim(G)/2}}e^{t\|\rho_\mathfrak{g}\|^2}\]
for $t\to0+$. This is Theorem~3 of Urakawa~\cite{urakawa}.

\item Theoretically each coefficient of the asymptotic expansion can be deduced using the expressions for the exact heat kernel available in the literatures (see for example \cite{urakawa}*{Thm.~2, p.~288}  or \cite{arede}*{Thm.~2.2, p.~57}). To my knowledge, the derivation of such expressions all rely on the structure theory of $G$. In contrast, our method illustrates that, as far as the asymptotic heat kernel is concerned, we need not know the structure theory, and we can get all the coefficients in a single setting.
\end{enumerate}
\end{rmk}

\begin{thm}\label{thm:asympkdhtk}
Let $r_t$ be the heat convolution kernel of $\KD(\mathfrak{g},\mathfrak{k})^2$. We have
\[ r_t^{\exp}\sim \frac{h_t^\mathfrak{g}}{j_\mathfrak{g}}e^{t\|\rho_{\mathfrak{k}}+\mu\|^2/2} \]
for $t\to 0+$, valid in a neighborhood of the origin.
\end{thm}
\begin{proof}
By Theorem~\ref{thm:genlapfml} and \ref{thm:duffltlap}, we have
\[ \KD(\mathfrak{g},\mathfrak{k})^2 = \frac12 \mathcal{Q}(\Delta_{\mathfrak{g}})  +\frac12\|\rho_{\mathfrak{k}}+\mu\|^2 \]
as differential operators on $(C^\infty(G)\otimes E)^K$. So 
\[ \exp(t\KD(\mathfrak{g},\mathfrak{k})^2) = \exp\bigl(t\mathcal{Q}(\Delta_{\mathfrak{g}})/2\bigr)\exp(t\|\rho_{\mathfrak{k}}+\mu\|^2/2).\]
Hence, 
\[ r_t= e^{t\|\rho_{\mathfrak{k}}+\mu\|^2/2}q_t, \]
where $q_t$ is the heat convolution kernel of $\mathcal{Q}( \Delta_{\mathfrak{g}})/2$. The asymptotic expansion for $r_t^{\exp}$ now follows from Proposition~\ref{prop:asymhtkduflap2}.
\end{proof}

\section{The Local Index Theorem}

We now come to the proof of the local index theorem for the cubic Dirac operator 
\[ \KD_{\mathfrak{g}/\mathfrak{k}}\ldef \eta^{-1}\circ\KD(\mathfrak{g},\mathfrak{k})\circ\eta\]
on $G/K$,  where $\eta\colon\Gamma(E(G))\to (C^\infty(G)\otimes E)^K$ is the linear isomorphism~\eqref{eq:equivfuncs}. As we have indicated in the introduction, the outline of our argument follows that of Berline and Vergne~\cite{berlinevergne}, but there are details that need to be verified.

We shall denote by $p_t$ the heat convolution kernel of $\KD_{\mathfrak{g}/\mathfrak{k}}^2$. It is possible to deduce the asymptotic expansion for $p_t$ from that of the heat convolution kernel $r_t$ of $\KD(\mathfrak{g},\mathfrak{k})^2$ by averaging over the fibers of $G\to G/K$; more precisely,
\begin{equation}
 p_t(e) = \int_K r_t(k)\nu(k)^{-1}\,dk.
 \label{eq:equivhtkcd}
\end{equation}
Equations like the above hold in general for Laplacians associated with principal bundles (see Berline et al.~\cite{bgv}*{Ch.~5} or Duistermaat~\cite{duistermaathk}*{Ch.~9}), and that is the main idea of Berline and Vergne's proof. The difference between their proof and ours, aside from the kind of Dirac operator used, is that whereas the principal bundle they use is the oriented orthonormal frame bundle, we use the principal bundle $G\to G/K$. Because of this, the calculations done in their proof do not quiet carry over. We shall refrain from presenting  calculations that are mere adaptations of theirs, but only those that needs to be verified independently in our case.

\begin{rmk}
\begin{enumerate}
\item Because $\eta$ and $\KD(\mathfrak{g},\mathfrak{k})$ are $G$-equivariant, $\KD_{\mathfrak{g}/\mathfrak{k}}$ is a $G$-equivariant differential operator on $\Gamma(E(G))$. 

\item There may be two pairs $(G,K)$ and $(H,L)$ such that $G/K$ is diffeomorphic to $H/L$. But if the Lie algebras of $G$ and $H$ are not isomorphic, then the corresponding cubic Dirac operators $\KD_{\mathfrak{g}/\mathfrak{k}}$ and $\KD_{\mathfrak{h}/\mathfrak{l}}$ may be different. An example (from Goette~\cite{goette}*{Rmk.~1.13, p.~8}) is when $(G,K)=(\Spin(4),\Spin(3))$ and $(H,L)=(\SU(2),\{e\})$.
\end{enumerate}
\end{rmk}

Changing the domain of the integral~\eqref{eq:equivhtkcd} via the exponential map of $K$ gives us
\[ p_t(e) = \int_{\mathfrak{k}} j_{\mathfrak{k}}^2(X)r_t^{\exp}(X)e^{-\nu_{*}(X)}\,dX.\]
Then, with our result for the asymptotic expansion of $r_t$ (Theorem~\ref{thm:asympkdhtk}), we have
\[p_t(e) \sim 
	e^{t\|\rho_{\mathfrak{k}}+\mu\|^2/2}\int_{\mathfrak{k}} \frac{j_{\mathfrak{k}}^2(X)}{j_{\mathfrak{g}}(X)}h_t^{\mathfrak{g}}(X) e^{-\nu_*(X)}
	\,dX.
\]
Now recall that the orthogonal decomposition $\mathfrak{g}=\mathfrak{k}\oplus \mathfrak{p}$ is $\ad(\mathfrak{k})$-invariant. For each $X\in\mathfrak{k}$, let $\ad^{\mathfrak{k}}(X)$ and $\ad^{\mathfrak{p}}(X)$ denote the restrictions of $\ad(X)$ to $\mathfrak{k}$ and $\mathfrak{p}$, respectively. Then the matrix of $\ad(X)=\ad^{\mathfrak{k}}(X)\oplus \ad^{\mathfrak{p}}(X)$ is block diagonal, so we have the factorization
\begin{equation*}
	j_{\mathfrak{g}}(X)= j_{\mathfrak{k}}(X)j_{\mathfrak{g}/\mathfrak{k}}(X), \label{eq:jkfactor}
\end{equation*}
where 
\begin{equation*} j_{\mathfrak{g}/\mathfrak{k}}(X) \ldef  \det\nolimits^{1/2}\biggl[ \frac{\sinh(\ad^{\mathfrak{p}}(X/2))}{\ad^{\mathfrak{p}}(X/2)}\biggr].
	\label{eq:hmspjacfctr}
\end{equation*}
Therefore,
\begin{equation}
  p_t(e) \sim 
	e^{t\|\rho_{\mathfrak{k}}+\mu\|^2/2} \int_{\mathfrak{k}} \frac{j_{\mathfrak{k}}(X)}{j_{\mathfrak{g}/\mathfrak{k}}(X)} h_t^{\mathfrak{g}}(X) e^{-\nu_*(X)}
	\,dX
\label{eq:htkrnlhmsid}
\end{equation}
for $t\to 0+$.

\begin{notation}
Let $\{Y_i\}_{i=1}^{\dim\mathfrak{p}}$ be an orthonormal basis for $\mathfrak{p}$. Define the linear map
\[ q\colon {\wedge(\mathfrak{p})}\to\Cl(\mathfrak{p}) \]
by
\[ q(Y_{i_1}\dotsb Y_{i_k}) = Y_{i_1}\dotsm Y_{i_k} \]
for any subset $\{Y_{i_1},\dotsc,Y_{i_k}\}$ of the basis. This is a vector space isomorphism, known as the \emph{Chevalley map}; it presents $\wedge(\mathfrak{p})$ as the associated graded algebra of the filtered algebra $\Cl(\mathfrak{p})$. We shall denote the composition of $\gamma^{\mathfrak{p}}\colon\mathfrak{k}\to\spin(\mathfrak{p})$ with $q^{-1}$  by
\[ \lambda^{\mathfrak{p}}\colon\mathfrak{k}\to \wedge(\mathfrak{p}).\]
That being so, applying $q^{-1}$ to Equation~\eqref{eq:decmopsospin2} gives us
\begin{equation}
\lambda^{\mathfrak{p}}(X) =-\frac{1}{2}\sum_{i,j=1}^{\dim\mathfrak{p}}\langle  X,[Y_i,Y_j]_{\mathfrak{g}}\rangle Y_iY_j\quad \in\wedge^2(\mathfrak{p}). \label{eq:adkpwg2}
\end{equation}
\end{notation}
\begin{rmk}\label{rmk:wg2cliffact}
We allow the image of $\lambda^{\mathfrak{p}}$ to act on the graded irreducible spinor space $S$ via the Chevalley map. This is useful when a product such as $\gamma^{\mathfrak{p}}(X_1)\dotsm\gamma^{\mathfrak{p}}(X_n)$ acts on $S$, and yet, we are only interested in the action ascribed to the part that is of top filtration order; in that case, we may replace $\gamma^{\mathfrak{p}}$ with $\lambda^{\mathfrak{p}}$ because 
\[ \gamma^{\mathfrak{p}}(X)\gamma^{\mathfrak{p}}(Y)=q(\lambda^{\mathfrak{p}}(X)\lambda^{\mathfrak{p}}(Y)) +(\text{terms of lower order}). \]
\end{rmk}

Another significance of $\lambda^{\mathfrak{p}}$ is that it is actually a curvature form in disguise. To reveal this hidden identity,  let $\Omega(G)$ denote the algebra of differential forms on $G$, and let $\Omega_\bas(G)$ be the subalgebra of basic forms, that is, the isomorphic image of $\Omega(G/K)$ under the pullback along the projection $\kappa\colon G\to G/K$. Recall that the curvature form $\Theta$ is an element of $\mathfrak{k}\otimes\Omega_\bas(G)$. And consider the following composition of isomorphisms:
\begin{equation}
 \Phi\colon \Omega_\bas(G) \xrightarrow{\ev_e}  {\wedge}(\mathfrak{p}^*) \xrightarrow{\sharp} \wedge(\mathfrak{p}),
\label{eq:basevm}
\end{equation}
where $\ev_e$ is the restriction to the tangent space at the identity, and $\sharp$ is the map induced by the inner product (so-called the ``raising of indices''). Then $\Phi(\Theta)$ is an element of $\mathfrak{k}\otimes\wedge(\mathfrak{p})$ that is essentially $\lambda^\mathfrak{p}$ as illustrated by the following lemma:
\begin{lem}\label{lem:lgrcurv}
Let $\theta$ be the canonical connection for the bundle $G\to G/K$ (see Remark on page~\pageref{rmk:equivweilact}). Let $\Theta$ be its curvature. Let $X$ be an arbitrary vector in $\mathfrak{k}$, and denote by $\Theta_X$ the basic form on $G$ obtained by taking the inner product of the $\mathfrak{k}$-factors of $\Theta$ with $X$. Then
\begin{equation}
 \lambda^\mathfrak{p}(X) = \Phi(\Theta_X), 
 \label{eq:lgrcurv}
\end{equation}
where $\Phi$ is the composite map~\eqref{eq:basevm}.
\end{lem}
\begin{rmk}
The liner map $\lambda^\mathfrak{p}\colon \mathfrak{k}\to\wedge(\mathfrak{p})$ can be identified as an element of $\mathfrak{k}^*\otimes\wedge(\mathfrak{p})$. Under the isomorphism $\mathfrak{k}^*\cong\mathfrak{k}$ induced by the inner product, $\lambda^\mathfrak{p}$ corresponds to the element
\begin{equation*}
 \Lambda = \sum_{i=1}^{\dim\mathfrak{k}}X_i\otimes\lambda^{\mathfrak{p}}(X_i),
\label{eq:lbardef}
\end{equation*}
where $\{X_i\}_{i=1}^{\dim\mathfrak{k}}$ is any orthonormal basis for $\mathfrak{k}$. Then Equation~\eqref{eq:lgrcurv} can be rewritten as
\begin{equation}
 \Lambda = \Phi(\Theta).
\label{eq:cubcurv}
\end{equation}
\end{rmk}\begin{proof}
Since $\Phi(\Theta_X)\in \wedge^2(\mathfrak{p})$, we can write 
\[ \Phi(\Theta_X)  = \frac12\sum_{i,j=1}^{\dim\mathfrak{p}} \langle X,\Theta(Y_i,Y_j)\rangle Y_i Y_j,\]
where $\{Y_i\}_{i=1}^{\dim\mathfrak{p}}$ is any orthonormal basis for $\mathfrak{p}$. By Cartan's structural equation $\Theta=d\theta+\frac12[\theta,\theta]_{\mathfrak{g}}$ (see Kobayashi and Nomizu~\cite{kobayashinomizu1}*{Thm.~5.2, p.~77}), we have:
\[
 \Theta(Y_i,Y_j) =d\theta(Y_i,Y_j)=-\theta([Y_i,Y_j]_{\mathfrak{g}})=-[Y_i,Y_j]_{\mathfrak{g}},
\]
where the last equality follows from the fact that $\theta$ is just the orthogonal projection of $\mathfrak{g}$ onto $\mathfrak{k}$ and that $[\mathfrak{p},\mathfrak{p}]\subset\mathfrak{k}$. Therefore,
\[
\langle X,\Theta (Y_i,Y_j)\rangle = -\langle   X,[Y_i,Y_j]_{\mathfrak{g}}\rangle,
\]
and hence,
\[ \Phi(\Theta_X)  = -\frac12\sum_{a,b=1}^{\dim\mathfrak{p}} \langle   X,[Y_i,Y_j]_{\mathfrak{g}}\rangle  Y_iY_j.\]
By comparing this expression with Equation~\eqref{eq:adkpwg2}, we obtain  Equation~\eqref{eq:lgrcurv}.
\end{proof}

Lemma~\ref{lem:jfuncwlgr} below captures how $\lambda^\mathfrak{p}$ (or its dual $\Lambda)$ enters into heat kernel calculations. Before stating the lemma, we make some preparatory remarks. Let $W$ be any finite-dimensional vector space over $\R$. We denote by $\wedge^+(W)$ the subalgebra of $\wedge(W)$ comprising all elements of even degree. As we identify $\mathfrak{k}$ with $\mathfrak{k}\otimes\{1\}\subset \mathfrak{k}\otimes\wedge^+(W)$,  any formal power series in $\mathfrak{k}^*$ extends uniquely to a map $\mathfrak{k}\otimes \wedge^+(W)\to \wedge^+(W)$. To wit, duality demands that the evaluation of any $\chi\in\mathfrak{k}^*$ against  an arbitrary element $\zeta=\sum X_i\otimes \zeta_i \in \mathfrak{k}\otimes\wedge^+(W)$ returns the value $\chi(\zeta)=\sum \chi(X_i)\zeta_i$; and the evaluation map $\ev_\zeta\colon \mathfrak{k}^*\to \wedge^+(W)$, $\chi\mapsto \chi(\zeta)$, extends uniquely as an algebra homomorphism to $\ev_\zeta\colon\R[[\mathfrak{k}^*]]\to \wedge^+(W)$; then $\ev_\zeta(\varphi)=:\varphi(\zeta)$ is the evaluation of the power series $\varphi$ at $\zeta$. All of this makes sense even when $W$ is replaced by $W_\C\ldef W\otimes\C$.

\begin{lem}\label{lem:jfuncwlgr} Let $h_t^{\mathfrak{k}}$ be the Gaussian on $\mathfrak{k}$. Let $\varphi$ be a smooth function on $\mathfrak{k}$ with sufficiently slow growth. Then the $\wedge(\mathfrak{p})$-valued function
\[ t  \mapsto   \int_{\mathfrak{k}} h_t^{\mathfrak{k}}(X) \varphi(X)e^{-\lambda^{\mathfrak{p}}(X)}\,dX\]
has an asymptotic expansion $\sum^\infty_{n=0}\Psi_nt^n$ for $t\to0+$. The $n$th coefficient $\Psi_n$ is contained in $\bigoplus_{q=0}^n\wedge^{2q}(\mathfrak{p})$. If $n\le \dim(\mathfrak{p})/2$, then the component of $\Psi_n$ of degree $2n$ (the highest degree part) is equal to that of the Taylor series of $\varphi$ evaluated at $-\Lambda$:
\begin{equation}  \Psi_{n}^{(2n)} = (\Taylor\varphi)(-\Lambda)^{(2n)}.
\label{eq:toppthmflshkex}
\end{equation}
\end{lem}
\begin{rmk}
This Lemma is similar in form to Lemma~11.3 in Duistermaat~\cite{duistermaathk}*{p.~137}. The proof given there can be carried over, but not quite. We will not give a separate proof but point out that in our case we need $\sum_{i=1}^{\dim\mathfrak{k}}\lambda^{\mathfrak{p}}(X_i)\lambda^{\mathfrak{p}}(X_i)=0$ for any orthonormal basis $\{X_{i}\}_{i=1}^{\dim \mathfrak{k}}$ for $\mathfrak{k}$. This equation can be verified using the Jacobi identity of the Lie bracket.
\end{rmk}

Equation~\eqref{eq:toppthmflshkex} is reminiscent of the construction of the Chern-Weil homomorphism, even more so in light of Equation~\eqref{eq:cubcurv}. The following proposition clearly exhibits their relationship. 

\begin{prop}\label{prop:smoothtospin}
Let $\theta$ denote the canonical connection for the principal $K$-bundle $\kappa\colon G\to G/K$. Let $\left.\Taylor\right|_\Lambda\colon C^\infty(\mathfrak{k})\to \wedge^+(\mathfrak{p}_\C)$ be the map that calculates the Taylor series of $\varphi\in C^\infty(\mathfrak{k})$ and then evaluates it at $\Lambda/2\pi i$. (We have inserted the extra factor $1/2\pi i$ to make certain formulas come out nicer afterward.) Let $\mathcal{A}$ denote the following composition of algebra homomorphisms:
\begin{equation*}
 C^\infty(\mathfrak{k})\xrightarrow{\left.\Taylor\right|_\Lambda} \wedge^+(\mathfrak{p}_\C) \xrightarrow{\Phi^{-1}}\Omega_\bas(G)_\C\xrightarrow{\kappa_*} \Omega(G/K)_\C,
 \label{eq:gencw}
\end{equation*}
where the homomorphisms $\Phi^{-1}$ and $\kappa_*$ are provided, respectively, by the inverse mappings of the composition~\eqref{eq:basevm} and the pullback isomorphism $\kappa^*\colon \Omega(G/K)\to \Omega_\bas(G)$. The restriction of $\mathcal{A}$ to $\mathfrak{k}$-invariant polynomials,
\[  S(\mathfrak{k}^*)^{\mathfrak{k}}\xrightarrow{\mathcal{A}}\Omega(G/K)_\C,\]
is precisely the Chern-Weil homomorphism induced by the canonical connection. 
\end{prop}
\begin{rmk}
The map $\left.\Taylor\right|_\Lambda$ factors through $S(\mathfrak{g}^*)$. 
\end{rmk}
\begin{proof}
Let us denote the Chern-Weil homomorphism as
\[ \CW\colon S(\mathfrak{k}^*)^{\mathfrak{k}}\to \Omega(G/K).\]
The assertion is that
\begin{equation*}
\operatorname{CW}(\varphi) =\kappa_{*}\Phi^{-1}\varphi(\Lambda/2\pi i)
\end{equation*}
for arbitrary $\varphi \in S(\mathfrak{k}^{*})^{\mathfrak{k}}$.
Since $\kappa_*$ and $\Phi^{-1}$ are isomorphisms, we may alternatively show that
\begin{equation}
 \varphi(\Lambda/2\pi i)= \Phi\,\kappa^*\!\CW(\varphi)
 \label{eq:cwa}
\end{equation}

First, we quickly recall how $\CW$ is constructed. Let $\Theta$ be the curvature of the canonical connection $\theta$. The curvature is a $\mathfrak{k}$-valued basic $2$-form on $G$, that is, a member of $\mathfrak{k}\otimes\Omega^+_\bas(G)$. The evaluation $\varphi(i\Theta/2\pi)$ of the polynomial $\varphi$ against $i\Theta/2\pi$ is in $\Omega_\bas^+(G)_\C$, so it is the pullback of a unique differential form in $\Omega(G/K)_\C$ along the canonical projection $G\to G/K$; that unique form on $G/K$ is by definition $\CW(\varphi)$. In short,
\[ \kappa^*\!\CW(\varphi)=\varphi(i\Theta/2\pi).\]
Therefore, 
\[ \Phi\kappa^*\!\CW(\varphi)=\Phi(\varphi(i\Theta/2\pi))=\varphi(i\Phi(\Theta)/2\pi),\]
where for the last equality we have used the fact that $\Phi$ is an algebra isomorphism. This proves, thanks to Equation~\eqref{eq:cubcurv}, our desired Equation~\eqref{eq:cwa}.
\end{proof}

\begin{prop}\label{prop:prinbexpjcb} Let $\mathcal{A}$ be as in Proposition~\ref{prop:smoothtospin}. Let $\varphi$ be the analytic function defined near the origin of $\mathfrak{k}$ by
\[  \varphi(X)= j_{\mathfrak{k}}(X)j_{\mathfrak{g}/\mathfrak{k}}^{-1}(X)\tr(e^{-\tau X}). \]
Then
\begin{equation} \mathcal{A}(\varphi)= \hat A (G/K;\theta) \ch(V(G);\theta), \label{eq:prinbexpjcb}
\end{equation} 
where $\hat A(G/K;\theta)$ and $\ch(V(G);\theta)$ are, respectively, the closed forms associated with $\theta$ (via the Chern-Weil homomorphism) in the Hirzebruch $\hat A$-class of $G/K$ and the Chern character of $V(G)$. If the rank of $K$ is strictly less than that of $G$, then $\hat A(G/K;\theta)=0$.  
\end{prop}
\begin{proof} 
Note that $\varphi$ is $\mathfrak{k}$-invariant, and so is its Taylor series. Now the Chern-Weil homomorphism $\CW$ extends to $\mathfrak{k}$-invariant  power series, but it still factors through $S(\mathfrak{k}^*)^\mathfrak{k}$; so by Proposition~\ref{prop:smoothtospin}, we have $\mathcal{A}=\CW$ on $\R[[\mathfrak{k}^*]]^\mathfrak{k}$. Hence,
\[ \mathcal{A}(\varphi)= \CW(j_{\mathfrak{k}})\CW(j^{-1}_{\mathfrak{g}/\mathfrak{k}})\CW(\tr(e^{-\tau})).\]
The differential forms $\CW(j_{\mathfrak{g}/\mathfrak{k}}^{-1})$ and $\CW(\tr(e^{-\tau}))$ are by definition $\hat A(G/K;\theta)$ and $\ch(V(G);\theta)$, respectively. As the last piece in establishing Equation~\eqref{eq:prinbexpjcb}, we claim that $\CW(j_{\mathfrak{k}})=1$. By definition, $\mathcal{A}=\kappa_*\circ\Phi^{-1}\circ\left.\Taylor\right|_\Lambda$; and  $\Taylor(j_\mathfrak{k})=j_\mathfrak{k}$; so it is sufficient to show that the evaluation of $j_\mathfrak{k}$ at $\Lambda/2\pi i$ is $1$. Now  
\[ j_{\mathfrak{k}}(\Lambda/2\pi i)=J(\ad^\mathfrak{k}(\Lambda)/2\pi i), \] 
where $J$ is the analytic function on $\so(\mathfrak{k})$ defined by $J(A) = \det^{1/2}(2\sinh(A/2)/A)$, and 
\[  \ad^\mathfrak{k}(\Lambda) \ldef  -\sum_{i=1}^{\dim\mathfrak{k}} \ad^\mathfrak{k}(X_i)\otimes \lambda^\mathfrak{p}(X_i)\in \so(\mathfrak{k})\otimes\wedge^+(\mathfrak{p}),\]
where  $\{X_i\}_{i=1}^{\dim\mathfrak{k}}$ is any basis for $\mathfrak{k}$ (the definition is independent of the choice). Now $\ad^\mathfrak{k}(\Lambda)$ is in fact zero; this is a consequence of the Jacobi identity of the Lie bracket. We omit the calculation (it is similar to that which we have remarked after Lemma~\ref{lem:jfuncwlgr}). From $\ad^\mathfrak{k}(\Lambda)=0$ it immediately follows that $J(\ad^\mathfrak{k}(\Lambda)/2\pi i)=1$, which proves our claim that $\CW(j_\mathfrak{k})=1$.

Finally, suppose $K$ is not of maximal rank in $G$, and let us see why $\hat A(G/K;\theta)$ vanishes. Let $J$ denote the same power series as before, but now viewed as a function on $\so(\mathfrak{p})$. Again, it is sufficient to show that $J(\ad^\mathfrak{p}(\Lambda)/2\pi i)=0$, or rather, that the function $J\circ\ad^\mathfrak{p}$ on $\mathfrak{k}$  is identically zero. And so we claim that when $K$ is not of maximal rank the representation $\ad^{\mathfrak{p}}\colon \mathfrak{k}\to \End(\mathfrak{p})$ is degenerate. To see why this is so, let $\mathfrak{t}$ be the Lie algebra of a maximal torus in $G$. We can choose the maximal torus so that it includes a maximal torus in $K$. Then the root spaces of $\mathfrak{k}$ are root spaces of $\mathfrak{g}$. If the rank of $K$ is strictly less than that of $G$, then $\mathfrak{s}\ldef \mathfrak{t}\cap\mathfrak{p}$ is nontrivial. But $[\mathfrak{k},\mathfrak{t}]\subset\mathfrak{k}$ and $[\mathfrak{k},\mathfrak{p}]\subset\mathfrak{p}$, so $[\mathfrak{k},\mathfrak{s}]\subset\mathfrak{k}\cap\mathfrak{p}=\{0\}$, which means that $\ad^{\mathfrak{p}}(\mathfrak{k})$ acts trivially on $\mathfrak{s}$. This proves that $\ad^\mathfrak{p}$ is degenerate, and we are done.
\end{proof}

All the necessary results to prove the local index theorem on $G/K$ are now at hand. Arranging them follows exactly as done in Berline and Vergne~\cite{berlinevergne}, so we only trace the outline in the proof.

\begin{thm}[The Local Index Theorem for $G/K$]\label{thm:hmglclindthm}
Let $G$ be a compact connected Lie group equipped with a bi-invariant metric. Let $K$ be a closed connected Lie subgroup of $G$. Let $S$ be the graded irreducible spinor space for $\Cl(\mathfrak{p})$, and let $\tau\colon \mathfrak{k}\to \End(V)$ be an irreducible representation of $\mathfrak{k}$ such that the $\mathfrak{k}$-action on $E\ldef S\otimes V$ via $\gamma^\mathfrak{p}\otimes \Id+\Id\otimes \tau$ lifts to a $K$-action. Consider the Kostant-Dirac operator $\KD_{\mathfrak{g}/\mathfrak{k}}$ on the bundle $E(G)\to G/K$. Assume that the bundle is $\Z/2\Z$-graded, relative to which $\KD_{\mathfrak{g}/\mathfrak{k}}$ is an odd operator. Let $p_t$ be the heat convolution kernel of $\KD_{\mathfrak{g}/\mathfrak{k}}^2$. 

If $K$ is of maximal rank in $G$, then the super trace of $p_t$, multiplied by the Riemannian volume form of $G/K$, satisfies
\[\Str(p_t)\,\vol  = \bigl.\hat A(G/K;\theta)\ch(V(G);\theta)\bigr|^{\topf} +O(t)
\]
for $t\to0+$, where the decoration $\left.\right|^{\topf}$ refers to the top degree part. If $K$ is not of maximal rank, then
\[ \Str(p_t)\,\vol  =  O(t).\]
\end{thm}
\begin{rmk}
It follows that  the graded index of $\KD_{\mathfrak{g}/\mathfrak{k}}$ is zero when $K$ is not of maximal rank in $G$. This verifies a result of Bott~\cite{bott}*{Thm.~II, p.~170} in our context.
\end{rmk}
\begin{proof}
Owing to the homogeneity of $G/K$, it is sufficient to prove for $\Str(p_t)$ at the identity coset $\bar e$. By Equation~\eqref{eq:indlierepclifm} and \eqref{eq:htkrnlhmsid},
\[ \Str( p_t(\bar e)) \sim
	e^{t\|\rho_{\mathfrak{k}}+\mu\|^2/2}\int_{\mathfrak{k}} h_t^{\mathfrak{g}}(X) j_{\mathfrak{k}}(X)j_{\mathfrak{g}/\mathfrak{k}}^{-1}(X)  \tr(e^{-\tau(X)})\Str(e^{-\gamma^{\mathfrak{p}}(X)})\,dX.  \]
The super trace of an element of $\Cl(\mathfrak{p})$ is nonzero only if the element is of top filtration order. So (as remarked on page~\pageref{rmk:wg2cliffact}) we may replace $\gamma^{\mathfrak{p}}$ with $\lambda^{\mathfrak{p}}$ in the integrand. We also note that $h_t^{\mathfrak{g}}(X)=h^{\mathfrak{k}}_t(X)/(2\pi t)^{d/2}$, where $d\ldef \dim(\mathfrak{p})$. Hence,
\begin{equation*}
 \Str(p_t(\bar e)) \sim \frac{e^{t\|\rho_{\mathfrak{k}}+\mu\|^2/2}}{(2\pi t)^{d/2}} \int_{\mathfrak{k}} h_t^{\mathfrak{k}}(X)  j_{\mathfrak{k}}(X)j_{\mathfrak{g}/\mathfrak{k}}^{-1}(X) \tr(e^{-\tau(X)}) \Str(e^{-\lambda^{\mathfrak{p}}(X)})\,dX. \label{eq:sptrkdirint}
\end{equation*}
Applying Lemma~\ref{lem:jfuncwlgr} with \[ \varphi(X) \ldef  j_{\mathfrak{k}}(X)j_{\mathfrak{g}/\mathfrak{k}}^{-1}(X)  \tr (e^{-\tau X}),\]
we conclude that
\begin{equation*}  \Str( p_t(\bar e)) \sim \frac{e^{t\|\rho_{\mathfrak{k}}+\mu\|^2/2}}{(2\pi t)^{d/2}}\sum_{n=0}^\infty \Str (\Psi_n)t^n, \label{eq:hmhtksptrasymp}
\end{equation*}
where $\Psi_n$ is contained in $\bigoplus_{q=0}^n\wedge^{2q}(\mathfrak{p})$. The super trace of an element in $\wedge^{2n}(\mathfrak{p})$ can be nonzero only if $2n$ is the top degree, that is, $2n=d$. Hence, if $d$ is odd, we have
\[ \Str(p_t(\bar e))\sim 0.\]
If $d$ is even (which includes the case where $K$ is of maximal rank), then 
\begin{equation*} 
\Str(p_t(\bar e)) =  \frac1{(2\pi)^{d/2}}\Str (\Psi_{d/2}) + O(t).	 \label{eq:straprxhtkrflhm}
\end{equation*}
It remains to check that 
\begin{equation}
 \frac1{(2\pi)^{d/2}}\Str (\Psi_{d/2})\,\vol = \bigl.\hat A(G/K;\theta)\ch(V(G);\theta)\bigr|^{\topf}. \label{eq:ldindnsty}
\end{equation}
Recalling Equation~\eqref{eq:toppthmflshkex}, we see that $\Psi_{d/2}^\topf=\varphi(\Lambda)^\topf=(2\pi i)^{d/2}\varphi(\Lambda/2\pi i)^\topf$. Consequently,
\[ 
\frac{1}{(2\pi)^{d/2}}\Str(\Psi_{d/2})\,\vol =\mathcal{A}(\varphi)^{\topf},
\]
where $\mathcal{A}$ is as in Proposition~\ref{prop:smoothtospin}. Our desired Equation~\eqref{eq:ldindnsty} now follows from  Proposition~\ref{prop:prinbexpjcb}.
\end{proof}

\section*{Acknowledgements}
This article is based on my doctoral dissertation. I am in debt of gratitude  to my advisor Nigel Higson. Thanks also to John Roe for the beneficial and instructive conversations. This research was partially supported under NSF grant DMS-1101382.

\end{document}